\documentclass[a4paper,11pt]{amsart}
\usepackage{amsmath,amsthm,amssymb,latexsym}

\newtheorem{thm}{Theorem}[section]
\newtheorem{lemma}[thm]{Lemma}

\newtheorem{rem}[thm]{Remark}

\newtheorem{prop}[thm]{Proposition}

\numberwithin{equation}{section}

\def\e{\varepsilon}

\def\L{\mathcal{L}}
\def\o{\omega}
\def\d{\delta}

\def\k{\kappa}

\def\p{\partial}
\def\O{\Omega}
\def\G{\Gamma}
\def\f{\varphi}
\def\R{\mathbb R}

\def\00{\infty}
\def\->{\rightarrow}

\input{amssym.def}

\title[two-phase flow problem]{Two-phase flow problem coupled \\with mean curvature flow}

\author[C. Liu]{Chun Liu}
\address[C. Liu]{Department of Mathematics, The Penn State
University, University Park, PA 16802 USA}
\email{liu@math.psu.edu}

\author[N. Sato]{Norifumi Sato}
\address[N. Sato]{Furano H.S., Furano, Hokkaido
076-0011 Japan}
\email{differ@mb.infosnow.ne.jp}

\author[Y. Tonegawa]{Yoshihiro Tonegawa}
\address[Y. Tonegawa]{Department of Mathematics, Hokkaido University, Sapporo 060-0810 Japan}
\email{tonegawa@math.sci.hokudai.ac.jp}

\date{}
\keywords{two-phase fluid, surface energy, varifold, phase field method}
\thanks{Y.T. is partially supported by JSPS Grant-in-aid for scientific research (B) 21340033}

\begin{document}

\begin{abstract}
We prove the existence of generalized solution for incompressible 
and viscous non-Newtonian two-phase fluid flow for spatial dimension 2 and 3. 
The phase boundary moves along with the fluid flow
plus its mean curvature while exerting surface tension force
to the fluid. An approximation scheme
combining the Galerkin method and the phase field method 
is adopted.
\end{abstract}

\maketitle

\makeatletter
\@addtoreset{equation}{section}%                     
\renewcommand{\theequation}{\thesection.\@arabic\c@equation}
\makeatother

\section{Introduction}
\quad 
In this paper we prove existence results for a problem on incompressible 
viscous two-phase fluid flow in the
torus $\O={\mathbb T}^d=({\mathbb R}/{\mathbb Z})^d$, 
$d=2,\,3$. A freely moving $(d-1)$-dimensional phase boundary
$\Gamma(t)$ separates the domain $\O$ into two domains 
$\O^+(t)$ and $\O^-(t)$, $t\geq 0$. The fluid flow is described 
by means of the velocity field $u:\O\times [0,\infty)\rightarrow
{\mathbb R}^d$ and the pressure 
$\Pi:\O\times [0,\infty)\rightarrow \R$. We assume the stress tensor
of the fluids is of the form $T^{\pm}(u,\Pi)=\tau^{\pm}(e(u))-\Pi\, I$ on $\O^{\pm}(t)$,
respectively. Here $e(u)$ is the symmetric part of the velocity gradient $\nabla u$, i.e. $e(u)=(\nabla u+\nabla u^T)/2$ and $I$ is the $d\times d$ identity matrix. 
Let $\mathbb{S}(d)$ be the set of $d\times d$ symmetric matrices. 
We assume that the functions $\tau^{\pm}:\mathbb{S}(d)\rightarrow\mathbb{S}(d)$ is locally Lipschitz and
satisfy for some $\nu_0>0$ and $p>\frac{d+2}{2}$ and for all $s,\,\hat{s}\in \mathbb{S}(d)$ 
\begin{equation}
\nu_0 |s|^p \leq \tau^{\pm}(s):s\leq \nu_0^{-1}(1+|s|^p),\label{taucond1}
\end{equation}
\begin{equation}
|\tau^{\pm}(s)|\leq \nu_0^{-1}(1+|s|^{p-1}),\label{taucond2}
\end{equation}
\begin{equation}
(\tau^{\pm}(s)-\tau^{\pm}(\hat{s})):(s-\hat{s})\geq 0.\label{taucond3}
\end{equation}
Here we define $A:B=\rm{tr}(AB)$ for $d\times d$ matrices $A,\, B$.
A typical example is $\tau^{\pm}(s)=(a^{\pm}+b^{\pm}|s|^2)^{\frac{p-2}{2}}s$
with $a^{\pm}>0$ and $b^{\pm}>0$.

We assume that the velocity field $u(x,t)$ satisfies the 
following non-Newtonian fluid flow equation:
\begin{eqnarray}
\frac{\p u}{\p t}+u\cdot\nabla u ={\rm div}\,(T^+(u,\Pi)),\hspace{.5cm}{\rm div}\, u=0 &\quad & {\rm on} \ \O^+(t), \ t> 0,\label{main1}\\
\frac{\p u}{\p t}+u\cdot\nabla u ={\rm div}\,(T^-(u,\Pi)),\hspace{.5cm}{\rm div}\, u=0 &\quad & {\rm on} \ \O^-(t), \ t> 0,\label{main2}\\
u^+= u^-,\hspace{.5cm}n\cdot (T^+(u,\Pi)-T^-(u,\Pi))=
\kappa_1 H &\quad & {\rm on} \ \Gamma(t), \ t> 0.\qquad \qquad \label{main3}
\end{eqnarray}
The upper script $\pm$ in \eqref{main3} indicates the limiting values 
approaching to $\Gamma(t)$ from $\O^{\pm}(t)$,
respectively, $n$ is the unit outer normal vector of $\p\O^+(t)$,
$H$ is the mean curvature vector of $\Gamma(t)$ and 
$\kappa_1>0$ is a constant. The condition \eqref{main3}
represents the force balance with an 
isotropic surface tension effect of the phase boundary. 
The boundary $\Gamma(t)$ is assumed to move with the velocity given by
\begin{equation}
V_{\Gamma}=(u\cdot n)n+\kappa_2 H \hspace{.5cm}{\rm on}
\quad\Gamma(t),\quad t> 0,
\label{velocity}
\end{equation}
where $\kappa_2>0$ is a constant. This differs from the
conventional kinematic condition ($\kappa_2=0$) and is 
motivated by the phase boundary motion with hydrodynamic
interaction. The reader is referred to \cite{Liu}
and the references therein for the relevant physical 
background. By setting $\f=1$ on $\O^+(t)$, $\f=-1$ on
$\O^-(t)$ and
\begin{equation*}
\tau(\f,e(u))=\frac{1+\f}{2}\tau^+(e(u))+\frac{1-\f}{2}
\tau^-(e(u))
\end{equation*}
on $\O^+(t)\cup\O^-(t)$, 
the equations \eqref{main1}-\eqref{main3} are expressed
in the distributional sense as
\begin{equation}
\begin{split}
\frac{\p u}{\p t}+u\cdot\nabla u &={\rm div}\,\tau(\f,e(u))
-\nabla \Pi +\kappa_1 H\mathcal{H}^{d-1}\lfloor_{\Gamma(t)} 
\hspace{.5cm} {\rm on} \ \O\times (0,\infty), \label{nsdist}\\
{\rm div}\, u&=0 \hspace{.5cm} {\rm on} \ \O\times (0,\infty).
\end{split}
\end{equation}
where $\mathcal{H}^{d-1}$ is the $(d-1)$-dimensional Hausdorff measure.
The expression \eqref{nsdist} makes it evident that the phase
boundary exerts surface tension force on the fluid wherever $H\neq 0$
on $\Gamma(t)$. Note that if $\Gamma(t)$ is a boundary of convex domain, the sign
of $H$ is taken so that the presence of surface tension tends to accelerate
the fluid flow inwards in general. 
We remark that the sufficiently smooth solutions of \eqref{main1}-\eqref{velocity} satisfy the following energy equality,
\begin{equation}
\frac{d}{dt}\left\{\frac{1}{2}\int_{\O}|u|^2\,dx+\kappa_1{\mathcal H}^{d-1}(\Gamma(t))\right\}=-\int_{\O}\tau(\f,e(u)):e(u)\,dx
-\kappa_1\kappa_2\int_{\Gamma(t)}|H|^2\,d{\mathcal H}^{d-1}. \label{energyeq}
\end{equation}
This follows from the first variation formula for the 
surface measure
\begin{equation}
\frac{d}{dt}{\mathcal H}^{d-1}(\Gamma(t))=-\int_{\Gamma(t)}
V_{\Gamma}\cdot H\, d{\mathcal H}^{d-1}
\label{firstvar}
\end{equation}
and by the equations \eqref{main1}-\eqref{velocity}. 

The aim of the present paper is to prove the time-global existence of
the weak solution for \eqref{main1}-\eqref{velocity} (see
Theorem \ref{maintheorem} for the precise statement). We construct
the approximate solution via the Galerkin method and the phase field method.
Note that it is not even clear for our problem 
if the phase boundary may stay as a 
codimension 1 object since a priori irregular flow field may tear
apart or crumble the phase boundary immediately, with a possibility
of developing singularities and fine-scale complexities. Even if we set 
the initial datum to be sufficiently regular, the eventual occurrence
of singularities of phase boundary or flow field may not be avoided
in general. To accommodate the presence of singularities of phase
boundary, we use the notion of varifolds from
geometric measure theory. In
establishing \eqref{velocity} we adopt the formulation due
to Brakke \cite{Brakke} where he proved the
existence of moving varifolds by mean curvature.
We have the extra transport effect $(u\cdot n)n$ which is 
not very regular in the present problem. Typically we would
only have $u\in L^p_{loc}([0,\infty);W^{1,p}(\Omega)^d)$. This poses a
serious difficulty in modifying Brakke's original construction
in \cite{Brakke} which is already intricate and involved. Instead
we take advantage of the recent progress on the understanding
on the Allen-Cahn equation with transport term to approximate the motion
law \eqref{velocity}, 
\[\frac{\p\f}{\p t}+u\cdot\nabla\f=\kappa_2\left(\Delta\f-\frac{W'(\f)}{\e^2}\right).\hspace{1cm}{\rm (ACT)} \]
Here $W$ is the equal depth double-well potential and
we set $W(\f)=(1-\f^2)^2/2$. When $\e\rightarrow 0$, we have proved
in \cite{LST1} that the interface moves according to the velocity \eqref{velocity}
in the sense of Brakke with a suitable regularity assumptions on $u$. 
To be more precise, we use a regularized
version of (ACT) as we present later for the result of \cite{LST1}
to be applicable. The result of \cite{LST1} was built upon those of
many earlier works, most relevant being \cite{Ilmanen1,Ilmanen2} which analyzed
(ACT) with $u=0$, and also \cite{Hutchinson,Tonegawa,Sato,Roeger}. 

Since the literature of two-phase flow is immense and continues to grow rapidly, we 
mention results which are closely related or whose aims point to some time-global existence
with general initial data. In the case without surface tension $(\kappa_1=\kappa_2=0)$, 
Solonnikov \cite{Solonnikov1} proved the time-local existence of classical solution. 
The time-local existence of weak solution was proved by Solonnikov \cite{Solonnikov2}, Beale \cite{Beale1}, Abels \cite{Abels1}, and others. For time-global existence of weak solution,
Beale \cite{Beale2} proved in the case that the initial data is small. Nouri-Poupaud \cite{Nouri} considered the case of multi-phase fluid. 
Giga-Takahashi \cite{GigaTakahashi} considered the problem within the framework of level set method. When $\kappa_1>0$, $\kappa_2=0$, Plotnikov \cite{Plotnikov} proved the 
time-global existence of varifold solution for $d=2$, $p>2$, and Abels \cite{Abels2} proved the time-global existence of measure-valued solution for 
$d=2, 3$, $p>\frac{2d}{d+2}$. When $\kappa_1>0$, $\kappa_2>0$, Maekawa \cite{Maekawa} proved the time-local existence of classical solution with $p=2$ (Navier-Stokes and Stokes)
and for all dimension. 
Abels-R\"{o}ger 
\cite{Abels-Roeger} considered a coupled problem of Navier-Stokes
and Mullins-Sekerka (instead of motion by mean curvature in the
present paper) and proved the existence of weak solutions. 
As for related phase field approximations of sharp interface model which we adopt in this paper, Liu and Walkington \cite{Liu} considered the case of fluids containing visco-hyperelastic particles. Perhaps the most closely related work to the
present paper is that of Mugnai and R\"{o}ger \cite{Mugnai} which studied 
the identical problem with $p=2$ (linear viscosity case) and $d=2,3$. There
they introduced the notion of $L^2$ velocity and showed that \eqref{velocity}
is satisfied in a weak sense different from that of Brakke for the limiting 
interface. 
Kim-Consiglieri-Rodrigues \cite{Kim} dealt with a coupling of Cahn-Hilliard and Navier-Stokes equations to describe the flow of non-Newtonian two-phase fluid 
with phase transitions. Soner \cite{Soner} dealt with a coupling of Allen-Cahn and heat equations to approximate the Mullins-Sekerka problem with kinetic undercooling.
Soner's work is closely related in that he showed the surface energy density bound
which is also essential in the present problem. 

The organization of this paper is as follows. In Section 2, we summarize the basic notations and main results. In Section 3 we construct a 
sequence of approximating solutions for the two-phase 
flow problem.
Section 4 describes the result of \cite{LST1}
which establishes the
upper density ratio bound for surface energy and which proves \eqref{velocity}. 
In the last Section 5 we combine the results from 
Section 3 and 4 and obtain the desired weak solution
for the two-phase flow problem. 

\section{Preliminaries and Main results}

\quad For $d\times d$ matrices $A,B$ we denote $A:B={\rm tr}\,(AB)$ and $|A|:=\sqrt{A:A}$. For $a \in \R^d$, 
we denote by $a\otimes a$ the $d\times d$ matrix with the $i$-th row and $j$-th
column entry equal to $a_i a_j$. 

\subsection{Function spaces}
\quad Set $\O={\mathbb T}^d$ throughout this paper. We set function spaces for $p>\frac{d+2}{2}$ as follows:
\begin{equation*}
\begin{split}
&{\mathcal V}=\left\{v \in C^{\infty}(\O)^d\,;\,{\rm div}\,v=0\right\},\\
&{\rm for} \ s\in {\mathbb Z}^+ \cup\{0\}, \ W^{s,p}(\O)=\{v \ : \ \nabla ^j
v\in L^p(\O) \ {\rm for } \ 0\leq j\leq s\},\\
&V^{s,p}= {\rm closure \ of} \ {\mathcal V} \ {\rm in \ the} \ 
W^{s,p}(\O)^d{\rm \mathchar`-norm.}
\end{split}
\end{equation*}
We denote the dual space of $V^{s,p}$ by $(V^{s,p})^*$. The $L^2$ inner
product is denoted by $(\cdot,\cdot)$.
Let $\chi_A$ be the characteristic function of $A$, and let $|\nabla\chi_A|$
be the total variation measure of the distributional derivative $\nabla \chi_A$.

\subsection{Varifold notations}
\quad We recall some notions from geometric measure theory and refer to \cite{Allard,Brakke,Simon} for more details. A {\it general $k$-varifold} in $\R^d$ is a 
Radon measure on $\R^d\times G(d,k)$, where $G(d,k)$ is the space of $k$-dimensional subspaces in $\R^d$. We denote the set of all general $k$-varifolds by ${\bf V}_k(\R^d)$. 
When $S$ is a $k$-dimensional subspace, we also use $S$ to denote the orthogonal projection matrix corresponding to
$\R^d\rightarrow S$. The first variation of $V$ can be written as
\begin{equation*}
\d V(g)=\int_{\R^d\times G(d,k)}\nabla g(x):S\,dV(x,S)
=-\int_{\R^d}g(x)\cdot H(x)\,d\|V\|(x) \quad {\rm if }\, \|\d V\|\ll \|V\|. 
\end{equation*}
Here $V \in {\bf V}_k(\R^d)$, $\|V\|$ is the mass measure of $V$, $g \in C_c^1(\R^d)^d$, $H=H_V$ is the generalized mean curvature vector if it exists and $\|\d V\|\ll \|V\|$
denotes that $\|\d V\|$ is absolutely continuous with respect to $\|V\|$. 

We call a Radon measure $\mu$ {\it $k$-integral} if $\mu$ is represented as $\mu=\theta{\mathcal H}^k\lfloor_X$, where $X$ is a countably $k$-rectifiable, ${\mathcal H}^k$-measurable set, 
and $\theta \in L^1_{\rm loc}({\mathcal H}^k\lfloor_X)$ is positive and 
integer-valued ${\mathcal H}^k$ a.e on $X$. ${\mathcal H}^k\lfloor_X$ denotes
the restriction of ${\mathcal H}^k$ to the set $X$. 
We denote the set of $k$-integral Radon measures 
by ${\mathcal{IM}}_k$. We say that a $k$-integral varifold is of {\it unit 
density} if $\theta=1$ ${\mathcal H}^k$ a.e. on $X$. For each such 
$k$-integral measure $\mu$ corresponds a unique $k$-varifold $V$
defined by 
\[\int_{\R^d\times G(d,k)}\phi(x,S)\,dV(x,S)=\int_{\R^d}\phi(x,T_x\mu)\,d\mu(x)\quad {\rm for} \ \phi\in C_c(\R^d\times G(d,k)),\]
where $T_x\mu$ is the approximate tangent $k$-plane. Note that $\mu=\|V\|$. We make such 
identification in the following. For this reason we define $H_{\mu}$ as $H_V$
(or simply $H$) if
the latter exists. When $X$ is a $C^2$ submanifold without 
boundary and $\theta$ is constant 
on $X$, $H$ corresponds to the usual
mean curvature vector for $X$. In the following we suitably
adopt the above notions on $\O={\mathbb T}^d$ such as
${\bf V}_k(\O)$, which present no essential difficulties.

\subsection{Weak formulation of free boundary motion}
For sufficiently smooth surface $\G(t)$ moving by the
velocity \eqref{velocity}, the following holds for
any $\phi\in C^2(\O;\R^+)$ due to the first variation formula
\eqref{firstvar}:
\begin{equation}
\frac{d}{dt}\int_{\G(t)}\phi\, d{\mathcal H}^{d-1}\leq
\int_{\G(t)}(-\phi H+\nabla\phi)\cdot\{\kappa_2 H+(u\cdot n)n\}\,
d{\mathcal H}^{d-1}.
\label{weakvelo}
\end{equation}
One can check that having this inequality for any $\phi\in C^2(\O;\R^+)$ implies
\eqref{velocity} thus \eqref{weakvelo} is equivalent to
\eqref{velocity}. Such use of non-negative test functions
to characterize the motion law is due to Brakke \cite{Brakke}
where he developed the theory of varifolds moving by the mean curvature.
Here we suitably modify Brakke's approach to incorporate the
transport term $u$. 
To do this we recall 
\begin{thm}{\bf (Meyers-Ziemer inequality)}
For any Radon measure $\mu$ on $\R^d$with 
\begin{equation*}D=\sup_{r>0,\, x \in {\mathbb R}^d}\frac{\mu(B_r(x))}{\o_{d-1}r^{d-1}}<\infty,
\end{equation*}
we have
\begin{equation}
\int_{\R^d}|\phi|\,d\mu\leq c D\int_{\R^d}|\nabla \phi|\,dx
\label{MZ1}
\end{equation}
for $\phi \in C_c^1(\R^d)$. Here $c$ depends only on $d$. \label{MZ}
\end{thm}
See \cite{Meyers} and \cite[p.266]{Ziemer}. 
By localizing \eqref{MZ1} to $\O={\mathbb T}^d$ we obtain (with $r$ in the
definition of $D$ above replaced by $0<r<1/2$) 
\begin{equation}
\int_{\O}|\phi|^2\, d\mu\leq c D \left(\|\phi\|_{L^2(\O)}^2+\|\nabla
\phi\|_{L^2(\O)}^2\right)
\label{MZ2}
\end{equation}
where the constant $c$ may be different due to the localization but 
depends only on $d$. 
The inequality \eqref{MZ2} allows us to define $\int_{\O}|\phi|^2\, d\mu$
for $\phi\in W^{1,2}(\O)$ by the standard density argument when $D<\infty$.

We define for any Radon measure $\mu$, $u\in L^2(\O)^d$ and $\phi\in C^1(\O:\R^+)$ 
\begin{equation}
{\mathcal B}(\mu,\, u,\, \phi)=\int_{\O}
(-\phi H+\nabla\phi)\cdot\{\kappa_2 H+(u\cdot n)n\}\, d\mu
\label{rhs}
\end{equation}
if $\mu\in {\mathcal{IM}}_{d-1}(\O)$ with generalized 
mean curvature $H\in L^2(\mu)^d$ and with
\begin{equation}\sup_{\frac12>r>0,\, x \in \O} \frac{\mu(B_r(x))}{\o_{d-1}r^{d-1}}<\infty
\label{den}
\end{equation} and $u\in W^{1,2}(\O)^d$. Due to the definition
of ${\mathcal {IM}}_{d-1}(\O)$, the unit normal vector $n$ is 
uniquely defined $\mu$ a.e. on $\O$ modulo $\pm$ sign. 
Since we have $(u,n)n$ in \eqref{rhs}, the choice of sign does not 
affect the definition.  
The right-hand side of \eqref{rhs} gives a well-defined 
finite value due to the stated conditions and \eqref{MZ2}. 
If any one of the conditions is not satisfied, we define ${\mathcal B}(\mu,\, u,\, \phi)=-\infty$. 

Next we note
\begin{prop}
For any $0<T<\infty$ and $p>\frac{d+2}{2}$,
$$\left\{u\in L^{p}([0,T];V^{1,p})\,;\,\frac{\p u}{\p t}\in L^{\frac{p}{p-1}}([0,T]; (V^{1,p})^*)\right\}\hookrightarrow C([0,T];\, V^{0,2}).$$
\label{embed}
\end{prop}
The Sobolev embedding gives $V^{1,p} \hookrightarrow V^{0,2}$ 
for such $p$ and we may apply \cite[p. 35, Lemma 2.45]{Malek} to obtain 
the above embedding. Indeed, we only need $p>\frac{2d}{d+2}$ for 
Proposition \ref{embed} to be and we have 
$\frac{d+2}{2}>\frac{2d}{d+2}$. Thus for this 
class of $u$ we may define $u(\cdot, t)\in V^{0,2}$ for all $t\in
[0,T]$ instead of a.e. $t$ and we may tacitly assume that we
redefine $u$ in this way for all $t$.
For $\{\mu_t\}_{t\in [0,\infty)}$, $u\in L^p_{loc}([0,\infty);V^{1,p})$ 
with $\frac{\p u}{\p t}\in L^{\frac{p}{p-1}}_{loc}
([0,\infty); (V^{1,p})^*)$ 
for $p>\frac{d+2}{2}$ and $\phi\in C^1(\O;\R^+)$, we define 
${\mathcal B}(\mu_t,\, u(\cdot,t),\, \phi)$ as in \eqref{rhs} for all $t\geq 0$.

\subsection{The main results}
Our main results are the following.
\begin{thm}
Let $d=2$ or $3$ and $p>\frac{d+2}{2}$. Let
$\O={\mathbb T}^d$. Assume that locally Lipschitz functions
$\tau^{\pm}:\mathbb{S}(d)\rightarrow\mathbb{S}(d)$ 
satisfy \eqref{taucond1}-\eqref{taucond3}. 
For any initial data $u_0\in V^{0,2}$ and 
$\O^+(0)\subset\O$ having $C^1$ boundary $\partial\O^+(0)$,
there exist
\begin{enumerate}
\item[(a)] $u \in L^{\infty}([0,\infty);V^{0,2})\cap L^p_{loc}([0,\infty);V^{1,p})$ 
with $\frac{\p u}{\p t}\in L^{\frac{p}{p-1}}_{loc}([0,\infty);(V^{1,p})^*)$, 
\item[(b)] a
family of Radon measures $\{\mu_t\}_{t\in [0,\infty)}$ with 
$\mu_t\in {\mathcal{IM}}_{d-1}$ for a.e. $t\in [0,\infty)$ and
\item[(c)] $\f \in BV_{loc}(\O\times [0,\infty)) \cap L^{\infty}([0,\infty);BV(\O)) \cap C^{\frac{1}{2}}_{loc}([0,\infty);L^1(\O))$ 
\end{enumerate}
such that the following properties hold:
\begin{enumerate}
\item[(i)] The triplet $(u(\cdot,t),\, \f(\cdot,t),\,\mu_t)_{t\in [0,\infty)}$ is a weak solution of \eqref{nsdist}. More precisely, for any $T>0$ we have
\begin{equation}
\int_0^T
\int_{\O}-u\cdot \frac{\p v}{\p t}+(u\cdot\nabla u)\cdot v+\tau(\f,e(u)):e(v)\,dxdt
=\int_{\O}u_0\cdot v(0)\,dx+\int_0^T\int_{\O}\kappa_1 H\cdot v \,
d\mu_t dt
\label{maintheorem1}
\end{equation} 
for any $v \in C^{\infty}([0,T];{\mathcal V})$ such that $v(T)=0$.
Here $H\in L^2([0,\infty);L^2(\mu_t)^d)$ is the generalized mean curvature vector 
corresponding to $\mu_t$. 
\item[(ii)] The triplet $(u(\cdot,t),\, \f(\cdot,t),\,\mu_t)_{t\in [0,\infty)}$
satisfies the energy inequality 
\begin{equation}
\begin{split}
\frac12\int_{\O}|u(\cdot,T)|^2\,dx+\kappa_1\mu_T(\O)&+\int_0^T\int_{\O}
\tau(\f,e(u)):e(u)\, dxdt+\kappa_1\kappa_2\int_0^T\int_{\O}|H|^2\, d\mu_t dt\\
&\leq \frac12\int_{\O}|u_0|^2\,dx+\kappa_1{\mathcal H}^{d-1}(\partial \O^+(0))
=: E_0\end{split}
\label{eneineq}
\end{equation}
for all $T<\infty$.
\item[(iii)] For all $0\leq t_1<t_2< \infty$ and $\phi\in C^2(\O;\R^+)$ we have
\begin{equation}
\mu_{t_2}(\phi)-\mu_{t_1}(\phi)\leq \int_{t_1}^{t_2}{\mathcal B}(\mu_t,\, u(\cdot,t),\, \phi)\, dt. 
\label{maintheorem3}
\end{equation}
Moreover, ${\mathcal B}(\mu_t,\, u(\cdot,t),\, \phi)\in L^{1}_{loc}([0,\infty))$.
\item[(iv)] We set $D_0=\sup_{0<r<1/2,\, x\in \O}\frac{{\mathcal H}^{d-1}
(\partial\O^+(0)\cap B_r(x))}{\o_{d-1}r^{d-1}}$. 
For any $0<T<\infty$, there exists a constant
$D=D(E_0,D_0,T,p,\nu_0,\k_1,\k_2)$ such that
$$\sup_{0<r<1/2,\,x\in \O}\frac{\mu_t(B_r(x))}{\omega_{d-1}r^{d-1}}
\leq D$$ for all $t\in [0,T]$.
\item[(v)] The function $\f$ satisfies the following properties.\\
 \ (1) $\f=\pm 1$ {\rm a.e. on} $\O$ for all $t\in 
[0,\infty)$.\\
 \ (2) $\f(x,0)=\chi_{\O^+(0)}-\chi_{\O\setminus\O^+(0)}$ {\rm a.e. on} $\O$.\\
 \ (3) ${\rm spt}|\nabla\chi_{\{\f(\cdot,t)=1\}}|
 \subset{\rm spt}\mu_t$ for all $t\in [0,\infty)$.
\item[(vi)] There exists 
\[T_1=T_1(E_0,D_0,p,\nu_0,\k_1,\k_2)>0\]
such that $\mu_t$ is of unit density for a.e. $t\in [0,T_1]$.
In addition $|\nabla\chi_{\{\f(\cdot,t)=1\}}|=\mu_t$ for a.e.
$t\in [0,T_1]$.
\end{enumerate} 
\label{maintheorem}
\end{thm}
%%%%%%%%%%%%%%%%%%%%%%%%%%%%%%%%%%%%%%%%
\begin{rem}
Somewhat different from $u=0$ case we do not expect that
\begin{equation}
\limsup_{\Delta t\rightarrow 0}\frac{\mu_{t+\Delta t}(\phi)
-\mu_t(\phi)}{\Delta t}\leq {\mathcal B}(\mu_t,\, u(\cdot,t),\phi)
\label{ve1}
\end{equation}
holds for all $t\geq 0$ and $\phi \in C^2(\O; \R^+)$ in
general. While
we know that the right-hand side is $<\infty$ (by definition) 
for all $t$, we
do not know in general if the left-hand side is $<\infty$.
One may even expect that at a time when 
$\int_{\O}|\nabla u(\cdot,t)|^p\,dx=\infty$, it may be $\infty$.
Thus we may need to define \eqref{velocity} in the integral form 
\eqref{maintheorem3}
for the definition of Brakke's flow. Note that in case $u=0$, 
one can show that the left-hand side of \eqref{ve1} is $<\infty$ for all $t\geq 0$
(see \cite{Brakke}).
\end{rem}
%%%%%%%%%%%%%%%%%%%%%%%%%%%%%%%%%%%%%%%%
\begin{rem}
The difficulty of multiplicities have been often encountered in the measure-theoretic
setting like ours. Varifold solutions constructed by Brakke \cite{Brakke} 
have the same properties in this regard. On the other hand, (vi) says that there is
no `folding' for some initial time interval $[0,T_1]$ at least. 
\end{rem}
\begin{rem}
In the following we set $\kappa_1=\kappa_2=1$ for notational simplicity, 
while all the argument can be modified with any positive $\kappa_1$ and 
$\kappa_2$ with no essential differences. On the other hand, their being
positive plays an essential role, and most of the estimates and claims deteriorate 
as $\kappa_1,\, \kappa_2\rightarrow 0$ and fail in the limit. How severely 
they fail in the limit may be of independent interest which we do not pursue
in the present paper. Note that $\kappa_2=0$ limit should correspond precisely 
to the setting of Plotnikov \cite{Plotnikov} for $d=2$.
\end{rem}

We use the following theorem. See \cite[p.196]{Malek} and the reference therein.

\begin{thm}{\bf(Korn's inequality)}
Let $1<p<\infty$. Then there exists a constant $c_K=c(p,d)$ such that
\[\|v\|_{W^{1,p}(\O)}^p\leq c_K (\|e(v)\|_{L^p(\O)}^p+\|v\|^p_{L^1(\O)})\]
holds for all $v \in W^{1,p}(\O)^d$.
\label{Korn}
\end{thm}

\section{Existence of approximate solution}
\quad In this section we construct a sequence of approximate solutions of
\eqref{main1}-\eqref{velocity} by 
the Galerkin method and the phase field method. The proof is a suitable modification 
of \cite{LinLiu} for the non-Newtonian setting even though we need to incorporate a
suitable smoothing of the interaction terms. 

First we prepare a few definitions. We fix a sequence $\{\e_i\}_{i=1}^{\infty}$ with $\lim_{i\rightarrow\infty}
\e_i=0$ and fix a radially symmetric non-negative function $\zeta\in C^{\infty}_c(\R^d)$ with ${\rm spt}\, \zeta\subset B_1(0)$ and
$\int\zeta\, dx=1$. For a fixed $0<\gamma<\frac12$ we define
\begin{equation}
\zeta^{\e_i}(x)=\frac{1}{\e_i^{\gamma}}\zeta\left(\frac{x}
{\e_i^{\gamma/d}}\right).
\label{zeta}
\end{equation}
We defined $\zeta^{\e_i}$ so that $\int \zeta^{\e_i}\, dx=1$,
$|\zeta^{\e_i}|\leq c(d)\e_i^{-\gamma}$ and $|\nabla\zeta^{\e_i}|
\leq c(d)\e_i^{-\gamma-\gamma/d}$.

For a given initial data $\O^+(0)\subset\O$ with $C^1$ 
boundary $\p \O^+(0)$, we can approximate $\O^+(0)$ in $C^1$ topology by a sequence
of domains $\O^{i+}(0)$ with $C^3$ boundaries. 
Let $d^{i}(x)$ be the signed distance
function to $\p \O^{i+}(0)$ so that $d^{i}(x)>0$ on $\O^{i+}(0)$ and
$d^{i}(x)<0$ on $\O^{i-}(0)$. Choose $b^{i}>0$ so that $d^{i}$ is $C^3$
function on the $b^{i}$-neighborhood of $\p\O^{i+}(0)$. 
Now we associate $\{\e_i\}_{i=1}^{\infty}$ with $\O^{i+}(0)$ by re-labeling 
the index if necessary so that $\lim_{i\rightarrow\infty}\e_i/b^i=0$ and
$\lim_{i\rightarrow\infty}\e_i^{j-1}|\nabla^j d^i|=0$ for
$j=2,\, 3$ on the $b^{i}$-neighborhood of $\p\O^{i+}(0)$. 
Let $h\in C^{\infty}(\R)$ be a function
such that $h$ is monotone increasing, $h(s)=s$ for $0\leq s\leq 1/4$
and $h(s)= 1/2$ for $1/2<s$, and define $h(-s)=-h(s)$ for $s<0$. 
Then define 
\begin{equation}
\f_0^{\e_i}(x)=\tanh(b^i h(d^i(x)/b^i)/\e_i).
\label{tanh}
\end{equation}
Note that we have $\f_0^{\e_i}\in
C^3(\O)$ and $\e_i^j|\nabla^j\f_0^{\e_i}|$ for $j=1,\, 2,\, 3$ are 
bounded uniformly independent of $i$. 
The well-known property of phase field approximation shows that 
\begin{equation}
\lim_{i\rightarrow
\infty}\|\f_0^{\e_i}-(\chi_{\O^+(0)}-\chi_{\O^-(0)})\|_{L^1(\O)}=0,\hspace{.5cm}
\frac{1}{\sigma}\left(\frac{\e_i|\nabla\f_0^{\e_i}|^2}{2}
+\frac{W(\f_0^{\e_i})}{\e_i}\right)\, dx\rightarrow {\mathcal H}^{d-1}\lfloor_{
\p\O^+(0)}
\label{tanhprop}
\end{equation}
as Radon measures. Here $\sigma=\int_{-1}^{+1}\sqrt{2W(s)}\, ds$. 

For $V^{s,2}$ with $s>\frac{d}{2}+1$ let $\{\o^i\}_{i=1}^{\infty}$ be a set of basis for $V^{s,2}$ such that it is orthonormal in $V^{0,2}$.
The choice of $s$ is made so that the Sobolev embedding theorem implies
$W^{s-1,2}(\O)\hookrightarrow L^{\infty}(\O)$ thus $\nabla \o^i \in L^{\infty}(\O)^{d^2}$.

Let $P_i:V^{0,2}\rightarrow V^{0,2}_i={\rm span}\,\{\o_1,\o_2,\cdots,\o_i\}$ be the orthogonal projection. 
We then project the problem \eqref{main1}-\eqref{velocity} to $V^{0,2}_i$ by utilizing the orthogonality in $V^{0,2}$.
Note that just as in \cite{LinLiu}, we approximate the mean
curvature term in \eqref{nsdist} by the appropriate
phase field approximation. We consider the following problem:
\begin{eqnarray}
\hspace{.3cm}
\frac{\p u^{\e_i}}{\p t}=P_i\left({\rm div}\,\tau(\f^{\e_i}, e(u^{\e_i}))-
u^{\e_i}\cdot\nabla u^{\e_i}-\frac{\e_i}{\sigma}{\rm div}\,((\nabla\f^{\e_i}\otimes\nabla\f^{\e_i})*\zeta^{\e_i})\right) & &  {\rm on} \ \O\times[0,\00),\label{appeq1}\\
u^{\e_i}(\cdot,t)\in V^{0,2}_i \qquad \qquad \qquad \qquad \qquad & & {\rm for} \ t\geq 0,\label{appeq2}\\
\frac{\p\f^{\e_i}}{\p t}+(u^{\e_i}*\zeta^{\e_i})\cdot\nabla\f^{\e_i}=\Delta\f^{\e_i}-\frac{W'(\f^{\e_i})}{\e_i^2} \qquad \qquad \quad & & {\rm on} \ \O\times [0,\00),\label{appeq3}\\
u^{\e_i}(x,0)=P_i u_0(x),\quad \f^{\e_i}(x,0)=\f_0^{\e_i}(x) \qquad \qquad \qquad & & {\rm on} \ \O.\label{appeq4}
\end{eqnarray}
Here $*$ is the usual convolution. 
We prove the following theorem.
\begin{thm}
For any $i\in {\mathbb N}$, $u_0\in V^{0,2}$
and $\f^{\e_i}_0$, there exists a weak solution $(u^{\e_i},\f^{\e_i})$ of 
\eqref{appeq1}-\eqref{appeq4} such that $u^{\e_i} \in L^{\infty}([0,\00);V^{0,2})\cap L^p_{loc}([0,\00);V^{1,p})$,
$|\f^{\e_i}|\leq 1$,   
$\f^{\e_i} \in L^{\infty}([0,\00);C^3(\O))$ and
$\frac{\p\f^{\e_i}}{\p t}\in L^{\infty}([0,\00);C^1(\O))$. \label{globalexistence}
\end{thm}
%%%%%%%%%%%%%%%%%%%%%%%%%%%%%%%%%%%%%%%%%%%%%
We write the above system in terms of $u^{\e_i}=\sum_{k=1}^{i}c^{\e_i}_k(t)\o_k(x)$ first. Since
\begin{gather*}
\left(\frac{d}{dt}u^{\e_i},\,\o_j\right)=\bigg(\frac{d}{dt}\sum_{k=1}^i c^{\e_i}_k(t)\,\o_k,\,\o_j\bigg)=\frac{d}{dt}c^{\e_i}_j(t),\\
(u^{\e_i}\cdot\nabla u^{\e_i},\,\o_j)=\sum_{k,l=1}^i c_k^{\e_i}(t)c_l^{\e_i}(t)(\o_k\cdot\nabla\o_l,\,\o_j),\\
\e_i({\rm div}\,((\nabla\f^{\e_i} \otimes\nabla\f^{\e_i})*\zeta^{\e_i}),\,\o_j)
= \,-\e_i \int_{\O} (\nabla\f^{\e_i}\otimes\nabla\f^{\e_i})*\zeta^{\e_i}:
\nabla \o_j \,dx,\\
\left({\rm div}\,\tau(\f^{\e_i},e(u^{\e_i})),\,\o_j\right)
= -\int_{\O}\tau(\f^{\e_i},e(u^{\e_i})):e(\o_j)\,dx
\end{gather*}
for $j=1,\cdots,i$, \eqref{appeq1} is equivalent to 
\begin{equation}
\begin{split}
\frac{d}{dt}c_j^{\e_i}(t)=& \,-\int_{\O}\tau(\f^{\e_i},e(u^{\e_i})):e(\o_j)\,dx
-\sum_{k,l=1}^i c_k^{\e_i}(t)c_l^{\e_i}(t)(\o_k\cdot\nabla\o_l,\,\o_j) \\ 
& +\frac{\e_i}{\sigma}\int_{\O}(\nabla\f^{\e_i}\otimes\nabla\f^{\e_i})*\zeta^{\e_i}:
\nabla \o_j \,dx= \,A^{\e_i}_j(t)+B_{klj} c^{\e_i}_k(t)c^{\e_i}_l(t)+D^{\e_i}_j(t).\label{appeq1-2}
\end{split}
\end{equation}
Moreover, the initial condition of $c_j^{\e_i}$ is 
\[c^{\e_i}_j(0)=(u_0,\,\o_j)\quad {\rm for} \ j=1,2,\dots,i.\]
We also set
\[E_0={\mathcal H}^{d-1}(\p\O^+(0))+\frac12 \int_{\O}|u_0|^2\, dx\]
and note that 
\begin{equation}
\frac{1}{\sigma}\int_{\O}\left(\frac{\e_i|\nabla\f_0^{\e_i}|^2}{2}
+\frac{W(\f^{\e_i}_0)}{\e_i}\right)\,dx+\frac12\sum_{j=1}^i(c^{\e_i}_j(0))^2\leq E_0
+o(1)
\label{eqeq}
\end{equation}
by \eqref{tanhprop} and by the projection $P_i$ being orthonormal. 

We use the following lemma to prove Theorem \ref{globalexistence}. 
\begin{lemma}
There exists a constant $T_0=T_0(E_0,i,\nu_0,p)>0$ such that \eqref{appeq1}-\eqref{appeq4} with \eqref{eqeq} has a weak solution $(u^{\e_i},\f^{\e_i})$ in $\O\times[0,T_0]$
such that $u^{\e_i} \in L^{\infty}([0,T_0];V^{0,2})\cap L^p([0,T_0];V^{1,p})$, $|\f^{\e_i}|\leq 1$, 
$\f^{\e_i} \in L^{\infty}([0,T_0];C^3(\O))$ and $\frac{\p\f^{\e_i}}{\p t}
\in L^{\infty}([0,T_0];C^1(\O))$. 
\label{localexistence}
\end{lemma}
{\it Proof.} Assume that we are given a function 
$u(x,t)=\sum_{j=1}^i c_j^{\e_i}(t)\o_j(x)\in C^{1/2}([0,T];V^{s,2})$ with
\begin{equation}
c^{\e_i}_j(0)=(u_0,\,\o_j),\hspace{.5cm}
\max_{t\in[0,T]}\left(\frac12\sum_{j=1}^i|c^{\e_i}_j(t)|^2\right)^{1/2}+
\sup_{0\leq t_1<t_2\leq T}\sum_{j=1}^i
\frac{|c_j^{\e_i}(t_1)-c_j^{\e_i}(t_2)|}{|t_1-t_2|^{1/2}}\leq \sqrt{2E_0}.\label{leraycond}
\end{equation}
We let $\f (x,t)$ be the solution of the following parabolic equation:
\begin{equation}
\begin{split}
\frac{\p\f}{\p t}+(u*\zeta^{\e_i})\cdot\nabla\f=\Delta\f-\frac{W'(\f)}{\e_i^2},\\
\f(x,0)=\f^{\e_i}_0(x).
\end{split}\label{acapprox}
\end{equation}
The existence of such $\f$ with $|\f|\leq 1$ is guaranteed by the standard theory of parabolic equations (\cite{Ladyzhenskaya}).  
By \eqref{acapprox} and the Cauchy-Schwarz inequality, we can estimate 
\begin{equation*}
\frac{d}{dt}\int_{\O}\left(\frac{\e_i|\nabla\f|^2}{2}+\frac{W(\f)}{\e_i}\right)\,dx 
\leq -\frac{\e_i}{2}\int_{\O} \left(\Delta\f-\frac{W'(\f)}{\e_i^2}\right)^2\,dx+\frac{\e_i}{2}\int_{\O} \left\{(u*\zeta^{\e_i})\cdot\nabla\f\right\}^2\,dx.
\end{equation*}
Since for any $t \in [0,T]$
\begin{equation*}
\|u*\zeta^{\e_i}\|^2_{L^{\infty}(\O)} \leq \e_i^{-2\gamma}\|u\|^2_{L^{\infty}(\O)} 
\leq i\e_i^{-2\gamma}\max_{1\leq j \leq i}\|\o_j(x)\|^2_{L^{\infty}(\O)}
\sum_{j=1}^i|c^{\e_i}_j(t)|^2 \leq c(i)E_0,
\end{equation*}
\begin{equation*}
\frac{d}{dt}\int_{\O}\left(\frac{\e_i|\nabla\f|^2}{2}+\frac{W(\f)}{\e_i}\right)\,dx
\leq c(i) E_0\int_{\O}\frac{\e_i|\nabla\f|^2}{2}\,dx.
\end{equation*}
This gives
\begin{equation}
\sup_{0\leq t \leq T}
\frac{1}{\sigma}\int_{\O}\left(\frac{\e_i|\nabla\f|^2}{2}+\frac{W(\f)}{\e_i}\right)\,dx \leq e^{c(i) E_0 T}E_0.\label{energyest}
\end{equation}
Hence as long as $T\leq 1$, 
\begin{equation}
|D_j^{\e_i}(t)| \leq c \|\nabla\o_j\|_{L^{\infty}(\O)}\frac{1}{\sigma}\int_{\O}\int_{\O}\e_i|\nabla\f(y)|^2\zeta^{\e_i}(x-y)\,dydx
\leq c(i)e^{c(i) E_0}E_0\label{Dest}
\end{equation}
by $\nabla\o_j \in L^{\infty}(\O)^{d^2}$ and  \eqref{energyest}.

Next we substitute the above solution $\f$ into the place of $\f^{\e_i}$, and solve \eqref{appeq1-2} with the initial condition $c^{\e_i}_j(0)=(u_0,\,\o_j)$. Since $\tau$
is locally Lipschitz with respect to $e(u)$, there is at least some short time $T_1$ such that \eqref{appeq1-2} has a unique solution $\tilde{c}^{\e_i}_j(t)$ on $[0,T_1]$ with the initial condition
$\tilde{c}_j^{\e_i}(0)=(u_0,\,\o_j)$ for $1\leq j\leq i$. We show that the solution exists up to $T_0=T_0(i,E_0,p,\nu_0)$ satisfying \eqref{leraycond}.
Let $\tilde{c}(t)=\frac12\sum_{j=1}^i|\tilde{c}^{\e_i}_j(t)|^2$.
Then,
\begin{equation*}
\frac{d}{dt}\tilde{c}(t)= 
A^{\e_i}_j\tilde{c}^{\e_i}_j+B_{klj}\tilde{c}^{\e_i}_k\tilde{c}^{\e_i}_l\tilde{c}^{\e_i}_j+D_j^{\e_i}\tilde{c}^{\e_i}_j.
\end{equation*}
By \eqref{taucond1} $A_j^{\e_i}\tilde{c}^{\e_i}_j\leq 0$ hence
\begin{equation*}
\frac{d}{dt}\tilde{c}(t)  \leq c(i,E_0)(\tilde{c}^{3/2}+\tilde{c}^{1/2}).
\end{equation*}
Therefore,
\begin{equation}
\arctan\sqrt{\tilde{c}(t)} \leq \arctan\sqrt{E_0}+2c(i,E_0) t.\label{arc}
\end{equation}
We can also estimate $|dc_j^{\e_i}/dt|$ due to \eqref{appeq1-2}, \eqref{Dest}, 
\eqref{arc} and
\eqref{taucond2} depending only on $E_0,i,p,\nu_0$. Thus, by choosing $T_0$ small depending only on $E_0,i,p,\nu_0$ we have the existence of solution
for $t\in[0,T_0]$ satisfying \eqref{leraycond}.
We then prove the existence of a weak solution on $\O\times [0,T_0]$ by using Leray-Schauder fixed point theorem (see \cite{Ladyzhenskaya}). We define
\[\tilde{u}(x,t)=\sum_{j=1}^i\tilde{c}^{\e_i}_j(t)\o_j(x)\]
and we define a map $\L:u\mapsto \tilde{u}$ as in the above procedure. Let
\begin{equation*}
\begin{split}V(T_0):=&\left\{u(x,t)
=\sum_{j=1}^i c_j(t)\o_j(x)\,;\,\,\max_{t\in[0,T_0]}\left(\frac12\sum_{j=1}^i|c_j(t)|^2\right)^{1/2}\right.\\ &\left.+
\sup_{0\leq t_1<t_2\leq T_0}\sum_{j=1}^i
\frac{|c_j(t_1)-c_j(t_2)|}{|t_1-t_2|^{1/2}}\leq \sqrt{2E_0},\,c_j(0)=(u_0,\,\o_j),\,c_j\in C^{1/2}([0,T_0]) \right\}.
\end{split}
\end{equation*}
Then $V(T_0)$ is a closed, convex subset of $C^{1/2}([0,T_0];V^{0,2}_i)$ equipped with the norm
\[\|u\|_{V(T_0)}=\max_{t\in[0,T_0]}\left(\frac12\sum_{j=1}^i|c_j(t)|^2\right)^{1/2}+
\sup_{0\leq t_1<t_2\leq T_0}\sum_{j=1}^i
\frac{|c_j(t_1)-c_j(t_2)|}{|t_1-t_2|^{1/2}}\]
and by the above argument $\L:V(T_0)\rightarrow V(T_0)$. Moreover by the Ascoli-Arzel\`a compactness theorem $\L$ is 
a compact operator. Therefore by using the Leray-Schauder fixed point theoremC$\L$ has a fixed point $u^{\e_i}\in V(T_0)$. We denote by $\f^{\e_i}$ the solution of \eqref{appeq3} and \eqref{appeq4}. Then $(u^{\e_i}, \f^{\e_i})$ is a weak solution of \eqref{appeq1}-\eqref{appeq4} in $\O\times [0,T_0]$. Note that we have the
required regularities for $\f^{\e_i}$ due to the regularity of 
$u^{\e_i}*\zeta^{\e_i}$ in $x$ and by the standard parabolic regularity theory.
$\hfill{\Box}$
%%%%%%%%%%%%%%%%%%%%%%%%%%%%%%%%%%%%%%%%%%%%%%%%%%%%%%
\begin{thm}
Let $(u^{\e_i},\f^{\e_i})$ be the weak solution of \eqref{appeq1}-\eqref{appeq4} 
with \eqref{eqeq} in $\O\times[0,T]$. Then the following energy estimate holds:
\begin{equation}
\begin{split}
\int_{\O}\frac{1}{\sigma}&\left(\frac{\e_i|\nabla\f^{\e_i}(\cdot,T)|^2}{2}+\frac{W(\f^{\e_i}(\cdot,T))}{\e_i}\right)+\frac{|u^{\e_i}(\cdot,T)|^2}{2}\,dx\\
&+\int_0^{T}\int_{\O}\frac{\e_i}{\sigma}\left(\Delta\f^{\e_i}-\frac{W'(\f^{\e_i})}{\e_i^2}\right)^2+\nu_0|e(u^{\e_i})|^p\,dxdt \leq E_0+o(1).
\label{localenergy1}
\end{split}
\end{equation}
Moreover for any $0\leq T_1<T_2<\infty$
\begin{equation}
\int_{T_1}^{T_2}\|u^{\e_i}(\cdot,t)\|_{W^{1,p}(\O)}^p\, dt\leq c_K
\{\nu_0^{-1}E_0+(T_2-T_1)E_0^{\frac{p}{2}}\}+o(1).
\label{localenergysup}
\end{equation}
\label{localenergy}
\end{thm}
%%%%%%%%%%%%%%%%%%%%%%%%%%%%%%%%%%%%%%%
{\it Proof.} 
Since $(u^{\e_i},\f^{\e_i})$ is the weak solution of \eqref{appeq1}-\eqref{appeq4}, we derive 
\begin{equation}
\begin{split}
& \frac{d}{dt}\int_{\O}\frac{1}{\sigma}\left(\frac{\e_i|\nabla\f^{\e_i}|^2}{2}+\frac{W(\f^{\e_i})}{\e_i}\right)+\frac{|u^{\e_i}|^2}{2}\,dx\\
& =\int_{\O}-\frac{\e_i}{\sigma}\frac{\p \f^{\e_i}}{\p t}\left(\Delta\f^{\e_i}-\frac{W'(\f^{\e_i})}{\e_i^2}\right)+\frac{\p u^{\e_i}}{\p t}\cdot u^{\e_i}\,dx\\
& =\int_{\O}-\frac{\e_i}{\sigma}\left(\Delta\f^{\e_i}-\frac{W'(\f^{\e_i})}{\e_i^2}-(u^{\e_i}*\zeta^{\e_i})\cdot\nabla\f^{\e_i}\right)
\left(\Delta\f^{\e_i}-\frac{W'(\f^{\e_i})}{\e^2}\right)\,dx\\
& +\int_{\O}\left\{{\rm div}\,\tau(\f^{\e_i},e(u^{\e_i}))-u^{\e_i}\cdot\nabla u^{\e_i} 
-\frac{\e_i}{\sigma}{\rm div}\,((\nabla\f^{\e_i}\otimes\nabla\f^{\e_i})*\zeta^{\e_i})\right\}\cdot u^{\e_i}\,dx=I_1+I_2.
\end{split}\label{localenergy1cal}
\end{equation}
Since ${\rm div}\, (u^{\e_i}*\zeta^{\e_i})=({\rm div}\, u^{\e_i})*\zeta^{\e_i}=0$,  
\begin{equation*}
\sigma I_1 = -\int_{\O}\e_i\left(\Delta\f^{\e_i}-\frac{W'(\f)}{\e_i^2}\right)^2\,dx+\e_i\int_{\O}(u^{\e_i}*\zeta^{\e_i})\cdot\nabla\f^{\e_i}\Delta\f^{\e_i}\,dx.
\end{equation*}
For $I_2$, with \eqref{taucond1} 
\begin{equation*}
\int_{\O}{\rm div}\,\tau(\f^{\e_i},e(u^{\e_i}))\cdot u^{\e_i}\,dx
=-\int_{\O}\tau(\f^{\e_i},e(u^{\e_i})):e(u^{\e_i})\,dx 
\leq -\nu_0\int_{\O}|e(u^{\e_i})|^p\,dx.
\end{equation*}
Moreover the second term of $I_2$ vanishes by ${\rm div}\,u^{\e_i}=0$ and
\begin{equation*}
\begin{split}
& 
-\int_{\O}\e_i {\rm div}\,(\nabla\f^{\e_i}\otimes\nabla\f^{\e_i}*\zeta^{\e_i})\cdot u^{\e_i}\,dx
= -\int_{\O}\e_i \left(\nabla \frac{|\nabla\f^{\e_i}|^2}{2}+
\nabla \f^{\e_i}\Delta\f^{\e_i}\right)*\zeta^{\e} \cdot u^{\e_i}\,dx\\
& = -\e_i\int_{\O}(u^{\e_i}*\zeta^{\e_i})\cdot\nabla\f^{\e_i}\Delta\f^{\e_i}\,dx.
\end{split}
\end{equation*}
Hence \eqref{localenergy1cal} becomes
\begin{equation*}
\frac{d}{dt}\int_{\O}\frac{1}{\sigma}\left(\frac{\e_i|\nabla\f^{\e_i}|^2}{2}+\frac{W(\f^{\e_i})}{\e_i}\right)+\frac{|u^{\e_i}|^2}{2}\,dx
\leq -\int_{\O}\frac{\e_i}{\sigma}\left(\Delta\f^{\e_i}-\frac{W'(\f^{\e_i})}{\e_i^2}\right)^2
+\nu_0 |e(u^{\e_i})|^p\, dx.
\end{equation*}
Integrating with respect to $t$ over $t\in[0,T]$ and by \eqref{eqeq}, we obtain \eqref{localenergy1}. The proof of \eqref{localenergysup} follows from \eqref{localenergy1} and Theorem \ref{Korn}.
$\hfill{\Box}$\\

{\it Proof of Theorem \ref{globalexistence}.}
For each fixed $i$ we have a short time existence for $[0,T_0]$ where $T_0$ depends only on $i,E_0,p,\nu_0$ at $t=0$. By Lemma \ref{localenergy} the energy at $t=T_0$ is again bounded by $E_0+o(1)$.
By repeatedly using Lemma \ref{localexistence}, Theorem \ref{globalexistence} follows.
$\hfill{\Box}$\\

\section{Proof of main theorem}

\quad In this section we first prove that $\{\f^{\e_i}\}_{i=1}^{\infty}$ in
Section 3 and the associated surface energy measures $\{\mu_t^{\e_i}\}_{
i=1}^{\infty}$ converge subsequentially to $\f$ and $\mu_t$ which satisfy
the properties described in Theorem \ref{maintheorem}. Most of the 
technical and essential ingredients have been proved in \cite{LST1}
and we only need to check the conditions to apply the results. We 
then prove that the limit velocity field satisfies the weak non-Newtonian
flow equation, concluding the proof of Theorem \ref{maintheorem}.

First we recall the upper density ratio bound of the surface energy.

\begin{thm} (\cite[Theorem 3.1]{LST1})
Suppose $d\geq 2$, $\O={\mathbb T}^d$, $p>\frac{d+2}{2}$, $\frac12>\gamma\geq 0$,
$1\geq \e>0$ and $\f$ satisfies
\begin{eqnarray}
\frac{\p \f}{\p t}+u\cdot\nabla\f=\Delta\f-\frac{W'(\f)}{\e^2} \qquad \qquad \quad & & {\rm on} \ \O\times [0,T],\label{allen1}\\
\f(x,0)=\f_0(x) \qquad \qquad \qquad & & {\rm on} \ \O,\label{allen2}
\end{eqnarray}
where $\nabla^i u,\, \nabla^j \f, \nabla^k \f_t\in C(\O\times[0,T])$ for $0\leq i,\, k\leq 1$ and $0\leq j\leq 3$. 
Let $\mu_t$ be the Radon measure on $\O$ defined by
\begin{equation}
\int_{\O}\phi(x)\, d\mu_t(x)=\frac{1}{\sigma}\int_{\O}\phi(x)\left(\frac{\e|\nabla\f(x,t)|^2}{2}+\frac{W(\f(x,t))}{\e}\right)\, dx
\label{dmu}
\end{equation}
for $\phi\in C(\O)$, where $\sigma=\int_{-1}^1
\sqrt{2 W(s)}\, ds$.
We assume also that
\begin{gather}
\sup_{\O}|\f_0|\leq 1\mbox{ and }\sup_{\O}\e^i|\nabla^i\f_0|\leq c_{1}\mbox{ for $1\leq
i\leq 3$},\label{inibound}\\
\sup_{\O}\left(\frac{\e|\nabla\f_0|^2}{2}-\frac{W(\f_0)}{\e}\right)\leq \e^{-\gamma},\label{disbd}\\
\sup_{\O\times[0,T]}\left\{\e^{\gamma}|u|,\, 
\e^{1+\gamma}|\nabla u|\right\}\leq c_{2}, \label{uinfbound}\\
\int_0^T\|u(\cdot,t)\|^p_{W^{1,p}(\O)}\, dt\leq c_3.\label{ubound}
\end{gather}
Define for $t\in [0,T]$
\begin{equation}
D(t)=\max\left\{\sup_{x\in\O,\, 0<r\leq \frac12}\frac{1}{\o_{d-1}r^{d-1}}
\mu_t(B_r(x)), 1\right\},\hspace{1.cm}D(0)\leq D_0.
\label{dtdef}
\end{equation}
Then there exist $\epsilon_1>0$ which depends only on
$d$, $p$, $W$, $c_1$, $c_2$, $c_3$, $D_0$, $\gamma$ and $T$,
and $c_4$ which depends only on $c_3$, $d$, $p$, $D_0$ and
$T$ such that for all $0<\e\leq \epsilon_1$, 
\begin{equation}
\sup_{0\leq t\leq T}D(t)\leq c_4.
\label{fin1}
\end{equation}
\label{mainmono}
\end{thm} 
Using this we prove
\begin{prop}
For $\{\f^{\e_i}\}_{i=1}^{\infty}$ in Theorem \ref{globalexistence}, define
$\mu_t^{\e_i}$ as in \eqref{dmu} replacing $\f$ by $\f^{\e_i}$, and define
$D^{\e_i}(t)$ as in \eqref{dtdef} replacing $\mu_t$ by $\mu_t^{\e_i}$. Given
$0<T<\infty$, there exists $c_5$ which depends only on $E_0,\, \nu_0, \,
\gamma,\, D_0,\, T,\, d,\, p$ and $W$ 
such that 
\begin{equation}
\sup_{0\leq t\leq T}D^{\e_i}(t)\leq c_5
\label{key}
\end{equation}
for all sufficiently large $i$. 
\label{du}
\end{prop}
{\bf Proof}. We only need to check the conditions of Theorem \ref{mainmono}
for $\f^{\e_i}$ and $\mu_t^{\e_i}$. Note that $u$ in \eqref{allen1} 
is replaced by $u^{\e_i}*\zeta^{\e_i}$. We have $d\geq 2$, $\O={\mathbb T}^d$,
$p>\frac{d+2}{2}$, $\frac12>\gamma\geq 0$, $1\geq\e>0$ and \eqref{allen1} and
\eqref{allen2}. The regularity of functions is guaranteed in Theorem 
\ref{globalexistence}. With an appropriate choice of $c_1$, \eqref{inibound}
is satisfied for all sufficiently large $i$ due to the choice of $\e_i$
in \eqref{tanh}. The sup bound \eqref{disbd} is satisfied with even 0 on 
the right-hand side instead of $\e_i^{-\gamma}$. The bound for $u^{\e_i}*
\zeta^{\e_i}$ \eqref{uinfbound} is satisfied due to \eqref{zeta} and
\eqref{localenergy1}, and \eqref{ubound} is satisfied due to 
\eqref{localenergysup}. Thus we have all the conditions, and Theorem \ref{mainmono}
proves the claim.
$\hfill{\Box}$

We next prove
\begin{prop}
For $\{u^{\e_i}*\zeta^{\e_i}\}_{i=1}^{\infty}$ in Theorem \ref{globalexistence}, there
exist a subsequence (denoted by the same index) and the limit $u\in
L^{\infty}([0,\infty);V^{0,2})\cap L^p_{loc}([0,\infty); V^{1,p})$ 
% with $\frac{\partial u}{\partial t}\in L^{\frac{p}{p-1}}_{loc}([0,\infty);
% (V^{1,p})^*)$ 
such that for any $0<T<\infty$
\begin{equation}
u^{\e_i}*\zeta^{\e_i}\rightharpoonup u\mbox{ weakly in }L^p([0,T]; W^{1,p}(\O)^d),
\hspace{1.cm}u^{\e_i}*\zeta^{\e_i}\rightarrow u\mbox{ strongly in }L^2([0,T];L^2(\O)^d).
\label{weak}
\end{equation}
\end{prop}
{\bf Proof}. 
Let $\psi \in V^{s,2}$ with $||\psi||_{V^{s,2}}\leq 1$. With \eqref{appeq1}, 
\eqref{appeq2} and integration by parts, we have
\begin{equation*}
\begin{split}
\left(\frac{\p u^{\e_i}}{\p t},\psi\right)&= 
\left(\frac{\p u^{\e_i}}{\p t}, P_i\psi\right) 
= \left(-u^{\e_i}\cdot \nabla u^{\e_i}+{\rm div}\,\tau(\f^{\e_i},e(u^{\e_i}))
-\frac{\e_i}{\sigma}{\rm div}
(\nabla\f^{\e_i}\otimes \nabla\f^{\e_i})*\zeta^{\e_i}, P_i\psi\right) \\
&=\left(u^{\e_i}\otimes u^{\e_i}-\tau(\f^{\e_i},e(u^{\e_i}))
+\frac{\e_i}{\sigma}
(\nabla\f^{\e_i}\otimes \nabla\f^{\e_i})*\zeta^{\e_i},\nabla P_i\psi\right).
\end{split}
\end{equation*}
Here we remark that 
\[\|\nabla P_i\psi\|_{L^{\infty}(\O)}\leq c(d) \|P_i\psi\|_{W^{s,2}(\O)}\leq c(d)\|\psi\|_{W^{s,2}(\O)}=c(d)\|\psi\|_{V^{s,2}}\leq c(d)\]
by $s> \frac{d+2}{2}$ and properties of $P_i$ (see \cite{Lions} or 
\cite[p.290]{Malek}).
Thus by \eqref{taucond2} and \eqref{localenergy1}, we obtain
\begin{equation*}
\left(\frac{\p u^{\e_i}}{\p t},\psi\right)\leq c(d,p,\nu_0)\left(1+E_0+\|u^{\e_i}
\|_{W^{1,p}(\O)}^{p-1}\right).
\end{equation*}
Again using \eqref{localenergy1} and integrating in time we obtain
\begin{equation}
\int_0^T\left|\left|\frac{\p u^{\e_i}}{\p t}\right|\right|_{(V^{s,2})^*}^{\frac{p}{p-1}}\,dt\leq c(d,p,E_0,\nu_0,T).
\label{utes}
\end{equation}
Now we use Aubin-Lions compactness Theorem \cite[p.57]{Lions} with 
$B_0=V^{s,2}$, $B=V^{0,2}\subset L^2(\O)^d$, $B_1=(V^{s,2})^*$, $p_0=p$ and $p_1=\frac{p}{p-1}$. 
Then there exists a subsequence still denoted by $\{u^{\e_i}\}_{i=1}^{\infty}$ such that 
\begin{equation*}
u^{\e_i} \rightarrow u \quad {\rm in} \ L^p([0,T];L^2(\O)^d). \label{u-converge1}
\end{equation*}
Since we have uniform $L^{\infty}([0,T];L^2(\O)^d)$ bound for $u^{\e_i}$, 
the strong convergence also holds in $L^2([0,T];L^2(\O)^d)$. 
Note that we also have proper norm bounds to extract weakly convergent 
subsequences due to \eqref{localenergy1}. For each $T_n$ 
which diverges to $\infty$ as
$n\rightarrow\infty$, we choose a subsequence and by choosing a 
diagonal subsequence, we obtain the convergent subsequence with \eqref{weak}
with $u^{\e_i}$ instead of $u^{\e_i}*\zeta^{\e_i}$. It is not difficult to
show at this point that the same convergence results hold for $u^{\e_i}*
\zeta^{\e_i}$ as in \eqref{weak}.
$\hfill\Box$

{\bf Proof of main theorem}.
At this point, the rest of the proof concerning the existence of 
the limit Radon measure $\mu_t$ and the limit $\f=\lim_{i\rightarrow
\infty}\f^{\e_i}$ and their respective properties described in Theorem
\ref{maintheorem} can be proved by almost line by line identical
argument in \cite[Section 4,\, 5]{LST1}.
The only difference is that the energy $E_0$ in \cite{LST1} depends also on
$T$, while in this paper $E_0$ depends only on the initial data due to
\eqref{localenergy1}. This allows us to have time-global
estimates such as $u\in L^{\infty}([0,\infty);V^{0,2})$ and
$\f\in L^{\infty}([0,\infty);BV(\O))$. The argument in \cite{LST1}
then complete the existence proof of Theorem \ref{maintheorem} 
(b), (c) along with (iii)-(vi). We still need to prove (a), (i) and (ii).  
 
Due to \eqref{utes}, \eqref{taucond2} and \eqref{localenergysup} we may extract a
further subsequence so that
\begin{equation}
\frac{\partial u^{\e_i}}{\partial t}
\rightharpoonup \frac{\partial u}{\partial t}\mbox{ weakly in }
L^{\frac{p}{p-1}}([0,T];(V^{s,2})^*),\hspace{.5cm}
\tau(\f^{\e_i},e(u^{\e_i}))\rightharpoonup \hat{\tau}
\mbox{ weakly in }L^{\frac{p}{p-1}}([0,T];L^{\frac{p}{p-1}}(\Omega)^{d^2}).
\label{meascon1}
\end{equation}
For $\omega_j\in V^{s,2}$ ($j=1,\cdots$) and $h \in C^{\infty}_c((0,T))$ we have
\begin{equation*}
\int_{\Omega}{\rm div}
((\nabla\f^{\e_i}\otimes\nabla\f^{\e_i})*\zeta^{\e_i})\cdot h\omega_j\,
dx=\int_{\Omega}\left(\Delta\f^{\e_i}-\frac{W'(\f^{\e_i})}{\e_i^2}\right)
\nabla\f^{\e_i}\cdot h\omega_j*\zeta^{\e_i}\, dx
\end{equation*}
by integration by parts and ${\rm div}\,\omega_j=0$. 
Thus the argument in \cite[p.212]{Lions} and the 
similar convergence argument in \cite{LST1} show
\begin{equation}
\int_0^T\left\{\left(\frac{\p u}{\p t},h\o_j\right)+\int_{\Omega}
(u\cdot\nabla u)\cdot h\omega_j+h\hat{\tau}:e(\omega_j)
\, dx\right\}dt=\int_0^T\int_{\Omega}H\cdot h\omega_j\, d\mu_t dt.
\label{meascon2}
\end{equation}
Again by the similar argument using the 
density ratio bound and Theorem
\ref{MZ} one show by the density argument and \eqref{meascon2} that
$\frac{\p u}{\p t}\in L^{\frac{p}{p-1}}([0,T];(V^{1,p})^*)$ and 
\begin{equation}
\int_0^T \left\{\left(\frac{\p u}{\p t}, v\right)+\int_{\Omega}(u\cdot\nabla u)\cdot v+\hat{\tau}:e(v)
\, dx\right\}dt=\int_0^T\int_{\Omega}H\cdot v\, d\mu_t dt.
\label{meascon3}
\end{equation}
for all $v\in L^p([0,T];V^{1,p})$. 
We next prove
\begin{equation}
\int_0^T\int_{\O}\hat{\tau}:e(v)\, dxdt
=\int_0^T\int_{\O}\tau(\f,e(u)):e(v)\, dxdt
\label{last}
\end{equation}
for all $v\in C^{\infty}_c((0,T);{\mathcal{V}})$. 
As in \cite[p.213 (5.43)]{Lions}, we may deduce that
\begin{equation}
\frac12 \|u(t_1)\|^2_{L^2(\Omega)}+\int_0^{t_1}\int_{\Omega}\hat{\tau}
:e(u)\, dxdt\geq \int_0^{t_1}\int_{\Omega}H\cdot u\, d\mu_t dt
+\frac12\|u(0)\|^2_{L^2(\Omega)}
\label{last1}
\end{equation}
for a.e. $t_1\in [0,T]$. 
We set for any $v\in V^{1,p}$ 
\begin{equation}
A_i^{t_1}=\int_0^{t_1}\int_{\Omega}
(\tau(\f^{\e_i},e(u^{\e_i}))-\tau(\f^{\e_i},e(v))):(e(u^{\e_i})-e(v))\, dxdt+
\frac12 \|u^{\e_i}(t_1)\|^2_{L^2(\Omega)}.
\label{last2}
\end{equation}
The property \eqref{taucond3} of $e(\cdot)$ shows that the first term of \eqref{last2}
is non-negative. We may further assume that $u^{\e_i}(t_1)$ 
converges weakly to $u(t_1)$ in $L^2(\Omega)^d$ thus we have
\begin{equation}
\liminf_{i\rightarrow\infty}A_i^{t_1}\geq\frac12 \|u(t_1)\|^2_{L^2(\Omega)}.
\label{last3}
\end{equation}
By \eqref{appeq1} we have
\begin{equation*}
\begin{split}
A_i^{t_1}=&\frac12\|u^{\e_i}(0)\|_{L^2(\Omega)}^2-\frac{\e_i}{\sigma}
\int_0^{t_1}\int_{\Omega}{\rm div}((\nabla\f^{\e_i}\otimes\nabla\f^{\e_i})
*\zeta^{\e_i})\cdot u^{\e_i}\\
&-\int_0^{t_1}\int_{\Omega}\tau(\f^{\e_i},e(u^{\e_i})):e(v)+
\tau(\f^{\e_i},e(v)):(e(u^{\e_i})-e(v))\, dxdt
\end{split}
\end{equation*}
which converges to
\begin{equation}
A^{t_1}=\frac12\|u(0)\|_{L^2(\Omega)}^2+\int_0^{t_1}\int_{\Omega}
H\cdot u\, d\mu_t dt-\int_0^{t_1}\int_{\Omega}\hat{\tau}:e(v)
+\tau(\f,e(v)):(e(u)-e(v))\, dxdt.
\label{last4}
\end{equation}
Here we used that $\f^{\e_i}$ converges to $\f$ a.e. on $\Omega\times
[0,T]$. By \eqref{last1}, \eqref{last3} and \eqref{last4}, we deduce that
\begin{equation*}
\int_0^{t_1}\int_{\Omega}(\hat{\tau}-\tau(\f,e(v))):
(e(u)-e(v))\, dxdt\geq 0.
\end{equation*}
By choosing $v=u+\epsilon\tilde{v}$, divide by $\epsilon$ 
and letting $\epsilon\rightarrow 0$, we prove \eqref{last}.
Finally, \eqref{eneineq} follows from \eqref{last}, the
strong $L^1(\O\times[0,T])$ convergence of $\f^{\e_i}$, 
the lower semicontinuity of the mean curvature square term
(see \cite{LST1}) and the energy equality appearing in Theorem \ref{localenergy}.
This concludes the proof of Theorem \ref{maintheorem}
$\hfill{\Box}$

\end{document}